\documentclass{article}

\usepackage{amssymb}
\usepackage{amsmath}
\usepackage{amscd}
\usepackage{amsthm}
\usepackage{indentfirst}

\title{Hopfian and co-Hopfian Modules over Artinian rings}
\author{F. C. Leary}
\date{}

\begin{document}
\newtheorem{thm}{Theorem}[section]
\newtheorem{cor}[thm]{Corollary}
\newtheorem{propn}[thm]{Proposition}
\newtheorem{axe}[thm]{Lemma}
\theoremstyle{definition}
\newtheorem{defn}[thm]{Definition}
\newtheorem{illus}[thm]{Example}
\newtheorem{remk}[thm]{Remark}

\maketitle

\begin{abstract}
If $R$ is a ring with $1$, we call a unital left $R$-module $M$  {\em Hopfian} ({\em co-Hopfian}) in the category of left $R$-modules if any epic (monic) $R$-module endomorphism of $M$ is an automorphism. In the case $R$ is a commutative Noetherian ring, we use a result of Matlis to characterize those injective $R$-modules that are co-Hopfian, and to characterize those that are Hopfian when $R$ is also reduced. We show that if $R$ is a commutative Artinian principal ideal ring, then an $R$ module $M$ is Hopfian (co-Hopfian) if and only if $M$ is finitely generated if and only if its injective envelope $E(M)$ is Hopfian (co-Hopfian) if and only if $E(M)$ is finitely generated. We note that the ``finite uniserial type" problem poses an obstacle to establishing this result for an arbitrary Artinian principal ideal ring.
\end{abstract}

The modules of the title are relatives of finite sets and finite dimensional vector spaces. If $R$ is a ring with $1,$ we call a unital left $R$-module $M$ {\em Hopfian}  in the category of left $R$-modules  if any surjective endomorphism of $M$ is an isomorphism (equivalently, if $M$ is not isomorphic to any of its proper quotients), and, dually, {\em co-Hopfian} if every injective endomorphism of $M$ is an isomorphism (equivalently, if $M$ is not isomorphic to any of its proper submodules). We will call $M$ {\em bi-Hopfian}  if it is both Hopfian and co-Hopfian (a workable name suggested by E. Wofsey). 

Any Hopf condition (Hopfian, co-Hopfian, or bi-Hopfian) can be considered as a type of finiteness condition because in the category of sets they are satisfied precisely by the finite sets. Clearly, a vector space over a field is bi-Hopfian if and only if it is finite dimensional.

According to R. Hirshon \cite{Hirshon}, quoting from \cite[p 415]{MKS}:``The question of whether a group is hopfian was first studied by Hopf, who using topological methods, showed that the fundamental groups of closed two-dimensional orientable surfaces are hopfian."  We believe the relevant paper of Hopf to be \cite{Hopf}.

Hiremath \cite{Hiremath} introduced the notion of Hopfian module and Hopfian ring. Soon, Varadarajan \cite{Var}, in a wide-ranging paper,  considered Hopfian and co-Hopfian objects among Boolean rings, function algebras, and compact manifolds. 

In \cite{Me}, we called co-Hopfian modules {\em Dedekind finite} ( and Hopfian modules {\em Dedekind cofinite}), not realizing that this was in conflict with standard usage (we were following  Stout \cite{stout}, who studied, in topoi, objects $X$ with the property that monomorphisms $X \rightarrow X$ were isomorphisms;  he chose the name Dedekind finite in keeping with Dedekind's definition of finite set in \cite{Ded}). In this paper, despite our firm belief that Dedekind finite/co-finite are the proper terms, we follow standard usage and call such modules  co-Hopfian and Hopfian.

In \cite{Wolmer2}, Vasconcelos  showed that if $R$ is a commutative ring, then any finitely generated $R$-module $M$ is  Hopfian (he credits Strooker \cite{Strooker} with establishing this result  independently). In particular, any commutative ring $R$ is Hopfian as a module over itself. Later \cite{Wolmer}, he proved that if $R$ is commutative, then every finitely generated $R$-module is co-Hopfian if and only if $R$ is $0$-dimensional. So, any Noetherian ring with this property must be Artinian. It seems natural to try to characterize those Artinian rings for which the co-Hopfian modules are precisely the finitely generated ones. We will determine a class of (commutative) Artinian rings for which this is true.

Our first section will set forth some notation, terminology, and facts that we will use throughout. Section two will use Matlis' structure theorem for injective modules over a commutative Noetherian ring to characterize the injectives with a Hopf condition (Hopfian, co-Hopfian, or bi-Hopfian).We also provide a characterization of Hopfian modules in the special cases in which the ring is an integral domain or is reduced. The third section shifts the emphasis to injective modules over a commutative Artinian ring. We see that the three Hopf conditions, along with several others, are equivalent for injective modules if the ring is semisimple. We continue by showing that these equivalent conditions hold for arbitrary modules over a commutative Artinian principal ideal ring, with the added feature that the conditions hold for a module if and only if they hold for its injective envelope.  We next examine what can be said if the Artinian ring is not necessarily commutative.

\section{Preliminaries}

We start by establishing our notation and terminology, most of which is fairly standard. We follow with some facts and theorems that will be useful to keep in mind. 

Throughout, $R$ is a ring with $1 \neq 0.$ All $R$-modules are unital. We use $E$ to denote an arbitrary  injective $R$-module, and $E(M)$ is the injective envelope of the $R$-module $M.$ We use $\oplus_nM$ to represent the direct sum of $n$ copies of the module $M,$ $n$ a cardinal number.

It is well-known that a finitely generated $R$-module over a Noetherian (resp. Artinian) ring is Hopfian (resp. co-Hopfian). More generally, any Noetherian (resp.Artinian) module is Hopfian (resp. co-Hopfian). We have the added feature that a finitely generated module over an Artinian ring is bi-Hopfian, since it has finite length. 

In \cite{Me}, we showed that for $R$ a commutative ring, the free module $R^n$ is co-Hopfian if and only if $n$ is a finite integer and $R$ is a {\em quoring} (regular elements are invertible), that is, $R$ is its own {\em total quotient ring}. We went on to show that $R^n$ is Hopfian if and only if $n$ is finite, a less restrictive requirement. So, for $R^n$ to be bi-Hopfian, we need both $R$ a quoring and $n$ finite. Compare with \cite[Prop 1.4, Thm3.2]{Var}.

Let $R$ be a commutative ring.  If only finitely generated $R$-modules are to be co-Hopfian, then $R$ must have only finitely many maximal ideals. For any maximal ideal $\mathfrak{m},$ the $R$-module $R/\mathfrak{m}$ is cyclic. If $\mathfrak{m} \neq \mathfrak{n}$ are maximal ideals of $R,$ then $\mathrm{Hom}(R/\mathfrak{m},R/\mathfrak{n})=0.$ Let $\Omega$ represent the maximal spectrum of $R.$ Then
$$M=\oplus_{\Omega}R/\mathfrak{m}$$
is co-Hopfian by part 2 of the following theorem. So, if $\Omega$ has infinite cardinality, then there is an $R$-module that is co-Hopfian but not finitely generated.

\begin{thm}[Thms 1.2-1.4 \cite{Me}]
Let $R$ be a ring.
\begin{enumerate}
\item Let $M$ be a Hofian (co-Hopfian) $R$-module. If $M$ decomposes as a direct sum of a family $\{M_i\}$ of nontrivial $R$-modules, then each $M_i$ is Hopfian (co-Hopfian).
\item Let $\{M_i\}$ be a family of nontrivial $R$-modules and suppose that $\mathrm{Hom}(M_i,M_j)=\mathrm{Hom}(M_j,M_i)=0$ whenever $i \neq j.$ If each $M_i$ is co-Hopfian, then so is $\oplus M_i.$ The result is valid with co-Hopfian replaced by Hopfian.
\item  Let $M$ be a Hofian (co-Hopfian) $R$-module. If $M$ decomposes as a direct sum of a family $\{M_i\}$ of nontrivial $R$-modules, then there are only finitely many summands isomorphic to a given $R$-module $N.$
\end{enumerate}
\end{thm}

\begin{remk}
We mentioned earlier that we were mistakenly referring to co-Hopfian modules as Dedekind finite in \cite{Me}. As a result, we were blissfully unaware of the fact that Hiremath \cite{Hiremath} had shown that summands of Hopfian modules are Hopfian (as in Theorem 1.1.1). We were also writing \cite{Me} contemporaneously with Varadarajan's paper \cite{Var}, and so duplicated some of his results.\qed
\end{remk}

Theorem 1.1.2  has a counterpart for indecomposable injectives, which can be deduced from \cite[Prop 4.21]{S-V}.

\begin{thm}
Let $R$ be Noetherian,  and $P,Q$ prime ideals of $R,$ neither of which contains the other. Then
$$\mathrm{Hom}(E(R/P),E(R/Q))=0.$$
\end{thm}

We can even say that $\mathrm{Hom}(R/P,R/Q)=0$ under the same hypotheses on $P$ and $Q.$ Note that $\mathrm{Hom}(R/P,R/Q)\cong (P:Q)/P.$ But in this instance, $(P:Q)=P$ (see \cite[Cor 3, p 11]{DGN}).

Certain changes of ring preserve both the Hopfian property and the co-Hopfian property. What we need is contained in the following theorem, pointed out to us by E. Enochs.

\begin{thm}
Let $R \rightarrow R^{\prime}$ be an epimorphism in the category of  rings. Let $M,N$ be $R^{\prime}$-modules. Then $\mathrm{Hom}_{R^{\prime}}(M,N)=\mathrm{Hom}_R(M,N).$
\end{thm}

This result is contained in Proposition XI.1.2 in \cite{Bo}. It can also be found  at the Stacks Project (stacks.math.columbia.edu/tag/04VM) in the form as written. We will use this theorem in the particular cases in which $R$ is commutative and $R^{\prime}$ is a quotient of $R,$ a localization of $R,$ or, if $R$ is a domain, the fraction field of $R.$

\section{The general injective module}

Let $R$ be a commutative Noetherian ring and $E$ an injective $R$-module. In his seminal paper \cite{Eben}, Matlis determined the structure of $E:$  for each prime ideal $P$ of $R$, there is a cardinal number $n_P$ such that $E \cong \bigoplus_P (\oplus_{n_P}E(R/P)).$ The $n_P$ characterize $E$ up to isomorphism. With this result in hand, we can obtain a characterization of the co-Hopfian injective $R$-modules. While this result does not allow us to classify the Hopfian injectives, we can do so in the special case in which $R$ is reduced, a case which includes integral domains, and so determine the bi-Hopfian modules over such $R.$

\begin{thm}
Let $E \cong \bigoplus_P (\oplus_{n_P} E(R/P))$ be an injective $R$-module. Then $E$ is co-Hopfian if and only if $n_P$ is finite for all $P.$
\end{thm}

\begin{proof}
The proof is straightforward.  By Theorem 1.1.3, $E$ co-Hopfian implies $n_P$ is finite for all $P.$  Conversely, let $f:E \rightarrow E$ be an injective $R$-module homomorphism. Then $f(E)$ is a nonzero injective summand of $E$ and is isomorphic to $E.$ If $f(E) \neq E,$ then there is a nonzero injective module $E^{\prime}$ such that $E = f(E) \bigoplus E^{\prime}.$ But $E^{\prime}$ is itself a direct sum of indecomposable injectives, so $E \cong \bigoplus_P (\oplus_{m_P}E(R/P))$ for cardinals $m_P$ with $m_P > n_P$ for at least one prime ideal $P.$ But, this is impossible by the result of Matlis cited previously. Thus, $f(E)=E$ and $E$ is co-Hopfian.
\end{proof}

We make three trivial observations. First, the class of co-Hopfian injective $R$-modules is closed under {\em injective} submodules. Secondly,  the class of co-Hopfian injective $R$-modules is closed under finite direct sums. Finally, an injective $R$-module that is not co-Hopfian has co-Hopfian injective submodules.

If $R$ is a field, its only prime ideal is the $0$-ideal and the theorem reduces to ``a vector space over a field is co-Hopfian if and only if it is finite dimensional." Of course, a vector space is Hopfian (and so bi-Hopfian) under the same conditions. Going forward, we will assume that $R$ is not a field.

For a general commutative Noetherian $R$ the theorem does not characterize Hopfian injectives. In order that $E$ be Hopfian, each summand $E(R/P)$ must be Hopfian by Theorem 1.1.1. But this is not always so, an example being  the divisible (hence, injective) abelian group $\mathbb{Z}(p^{\infty}),$ the Pr\"{u}fer $p$-group. So, we should decide for which primes $P$ the injective module $E(R/P)$ is Hopfian. The notion of divisible module will assist  us. The following definition of divisible module is more or less standard (see \cite[Def 3.16]{Lam2} for an alternative definition that has its proponents).
 
\begin{defn}
Let $R$ be a commutative ring and $M$ an $R$-module. We say $M$ is {\em divisible} if $rM=M$ for each non-zerodivisor of $R.$
\end{defn}

It is well-known that over any commutative ring $R,$ injective modules are divisible. We are already assuming that $R$ is Noetherian. If $R$ is also a domain, then we see that $E(R/P)$ is not Hopfian for any nonzero prime $P.$ The key is that the nonzero elements of any nonzero prime $P$ are not zerodivisors of $R.$ Hence, $rE(R/P)=E(R/P)$ for any nonzero $r \in R.$ But $P$ is the unique associated prime of $E(R/P).$ So $P$ is the annihilator of some nonzero $x \in E(R/P).$ Hence, for any  nonzero $r \in P,$ multiplication by $r$ provides a surjective homomorphism of $E(R/P)$ that is not an isomorphism, and so $E(R/P)$ is not Hopfian. Of course, if $P$ is the $0$-ideal, the minimal prime ideal of $R,$ then $E(R/P)=E(R)=K,$ the fraction field of $R,$ which is Hopfian. The following theorem summarizes.

\begin{thm}
Let $R$ be a commutative Noetherian domain. The only Hopfian injective $R$-modules are the modules $K^n,$ where $K$ is the fraction field of $R$ and $n$ is a positive integer. These $K^n$ are also the only bi-Hopfian injective $R$-modules.
\end{thm}

\begin{proof}
Since $K^n$ is a $K$-vector space, it is bi-Hopfian as a $K$-module. But then it is also bi-Hopfian as an $R$-module (Theorem 1.4).
\end{proof}

Returning to commutative Noetherian $R$ we have the following result, which reduces to Theorem 2.3 when $R$ is a domain.

\begin{thm} 
If $P$ be a minimal prime of $R,$ then $E=\oplus_n E(R/P)$ is bi-Hopfian for any positive integer $n.$
\end{thm}

We know $E$ is co-Hopfian by Theorem 2.1, but we can prove $E$ is bi-Hopfian directly.

\begin{proof}
$E(R/P)$ is also an $R_P$-module, and, as such, is isomorphic to the $R_P$-module $E(R_P/PR_P).$ But $R_P$ is local Artinian with maximal ideal $PR_P,$ so $E(R_P/PR_P)$ is finitely generated \cite[Cor 3.86]{Lam2} as an $R_P$-module. Hence  $E(R_P/PR_P)$ has finite length, and so is both Hopfian and co-Hopfian as an $R_P$-module and so is bi-Hopfian as an $R$-module (Theorem 1.4). Now note that $\oplus_n E(R/P)$ also has finite length as an $R_P$-module, so that it is bi-Hopfian as an $R_P$-module, and so bi-Hopfian as an $R$-module.
\end{proof}

\begin{cor}
Let $E=\oplus(\oplus_{n_P}E(R/P)),$ the exterior sum over the minimal primes of $R$ and the $n_P$ nonnegative integers. Then $E$ is bi-Hopfian.
\end{cor}

\begin{proof}
If $E$ is to be bi-Hopfian, then the $n_P$ must be finite by Theorem 1.1.3. Hence, each $\oplus_{n_P}E(R/P)$ is bi-Hopfian by the theorem. If $P,Q$ are distinct minimal primes, then Theorem 1.3 yields that both $\mathrm{Hom}(E(R/P),E(R/Q))$ and $\mathrm{Hom}(E(R/Q),E(R/P))$ are $0,$ so that $E$ is bi-Hopfian by Theorem 1.1.2.
\end{proof}

We can actually say a bit more. If $R$ is not a domain, then the zero ideal is not prime. But, since $R$ is Noetherian, $0$ has a primary decomposition, and the elements of $R$ which are not zerodivisors are those in the complement of the union of the finitely many associated primes of $0.$ Call these primes $P_1, \ldots , P_n.$ If $Q$ is a prime distinct from the $P_i,$ then prime avoidance gives us an element $q \in Q$ that does not lie in any $P_i.$ Hence, $q$ is not a zero divisor of $R,$ and $qE(R/Q)=E(R/Q)$ since $E(R/Q)$ is divisible. But, again, $q$ annihilates some nonzero element of $E(R/Q), $ so multiplication by $q$ provides a surjective endomorphism of $E(R/Q)$ which is not an isomorphism. Hence, $E(R/Q)$ is not Hopfian. 

Now, if $R$ is reduced (that is, $R$ has no nonzero nilpotents), then the associated primes of $0$ are the minimal primes of $R,$ and the set of zerodivisors of $R$ is the union of them. We now have a clear dichotomy: if $P$ is a minimal prime, $E(R/P)$ is bi-Hopfian; for other primes $P,$ $E(R/P)$ is not bi-Hopfian. Arguing as in the  proof of Corollary 2.5, we see that Proposition 2.4 has an additional corollary.

\begin{cor}
Let $R$ be a reduced Noetherian ring. An injective $R$-module $E$ is bi-Hopfian if and only if $E=\oplus(\oplus_{n_P}E(R/P)),$ the outer sum taken over the minimal primes of $R,$ with the $n_P$ non-negative integers.
\end{cor}

\begin{remk}
Note that all the results of this section impose finiteness conditions on an injective module $E.$ The cardinals $n_P$ determining the number of summands $E(R/P)$ must all be finite. This situation will continue to occur in the sequel.
\end{remk}

\section{Application to Artinian rings}

Throughout this section, $R$ is a commutative Artinian ring. For such an $R,$ the maximal primes are also minimal primes, and the situation improves dramatically: we get a characterization of  Hopfian injectives as in Theorem 1.1, and more besides.  We discover that an injective $R$-module $E$ is Hopfian if and only if it is co-Hopfian if and only if it is finitely generated and other equivalent conditions as well. For a commutative Artinian ring $R$, indecomposable injectives are finitely generated (this is no longer the case if $R$ is not commutative \cite[Example 2.3]{Krause} \cite[p 375]{R-Z}). If $R$ is also a principal ideal ring ({\em PIR}), the equivalences mentioned earlier hold for  $R$-modules $M,$ even if not injective. In this case, there is also a close bond between an $R$-module $M$ and its injective envelope $E(M):$ $M$ is Hopfian if and only if $E(M)$ is.

\begin{thm}
Let $E$ be an injective $R$-module.  Then $E$ is Hopfian if and only if it is co-Hopfian.
\end{thm}

\begin{proof}
 Let $E=\oplus_P(\oplus_{n_p}E(R/P))$ be an injective $R$-module.

($\Leftarrow$) Let each $n_P$ be finite. Since $R$ has only finitely many prime ideals, all maximal, there are only finitely many summands, and each is finitely generated \cite[Cor 3.86]{Lam2}. Thus $E$ is finitely generated, and so has finite length (since $R$ is also Noetherian). Hence, $E$ is Hopfian. 

($\Rightarrow$) By Theorem 1.1, each $E(R/P)$ must be Hopfian, which it is by Theorem 2.4, and appear only finitely many times. Hence, each $n_P$ is finite. Hence, $E$ is co-Hopfian by Theorem 2.1.
\end{proof}

So for an injective $R$-module $E,$ Hopfian, co-Hopfian, and bi-Hopfian are equivalent conditions. We have an immediate corollary that follows from the proof of the theorem.

\begin{cor}
An injective $R$-module $E$ is Hopfian (hence co-Hopfian) if and only if $E$ is finitely generated.
\end{cor}

The following result, which is essentially \cite[Thm 3.64(3)]{Lam2} introduces a link between module and injective envelope.

\begin{thm}
Let $M$ be an $R$-module. Then $M$ is finitely generated if and only if $E(M)$ is finitely generated.
\end{thm}

\begin{proof}
($\Rightarrow$) If $M$ is finitely generated, then $E(M)$ is the direct sum of finitely many indecomposable injectives $E(R/P),$ $P$ a prime of $R,$ \cite[Thm 3.48]{Lam2}\cite[Thm 4.7]{S-V}, and these $E(R/P)$ are finitely generated. So, $E(M)$ is finitely generated. 
($\Leftarrow$) Clear since $R$ is Artinian.
\end{proof}

 Before proceeding to our next result, we introduce the concept of  finitely co-generated module. There are several definitions possible. See \cite[Prop 19.1]{Lam2} for four equivalent versions. For our purposes, the definition of Vamos, which appears in \cite[Sec 3.4]{S-V} with the name {\em finitely embedded}, is appropriate.

\begin{defn}
An $R$-module $M$ is {\em finitely co-generated} if its injective envelope $E(M)$ is isomorphic to the direct sum of the injective envelopes of finitely many simple $R$-modules.
\end{defn}

The theorem that follows is essentially an extension of the Hopkins-Levitzki Theorem to injective modules over an Artinian ring (\cite[Thm 15.20, Cor 15.21]{A-F}\cite[Thm 4.15]{Lam1}).

\begin{thm}
Let $E$ be an injective $R$-module. The following are equivalent.
\begin{enumerate}
\item $E$ is Hopfian.
\item $E$ is co-Hopfian.
\item $E$ is finitely generated.
\item $E$ has finite length (has a finite composition series).
\item $E$ has ACC (is Noetherian).
\item $E$ has DCC (is Artinian).
\item $E$ is finitely cogenerated.
\end{enumerate}
\end{thm}

\begin{proof}
The equivalence of $(3)-(6)$ follows from Hopkins-Levitzki (for {\em any} $R$-module). For the equivalence of (6) and (7) see \cite[Prop 19.4]{Lam2} or \cite[Thm 3.21]{S-V}. The equivalences of (1)-(3) is Corollary 3.2.
\end{proof}

\begin{remk}
An additional condition: (8) $E$ is bi-Hopfian, is implicit in the theorem. So, any one of the Hopf conditions (Hopfian, co-Hopfian, bi-Hopfian) is equivalent to the others. 

If an injective $E$ has a Hopf condition, then any submodule $M$ of $E$ will be finitely generated, and so have finite length, since $R$ is Artinian. Hence, all submodules of $E$ have that Hopf condition. In the terminology of Vedadi \cite[p 144]{Vedadi}, $E$ is {\em totally Hopfian} (or {\em totally co-Hopfian}, or {\em totally bi-Hopfian}).

Any finitely generated $R$-module is co-Hopfian, but the converse is {\em not} true, as we shall see in Example 3.15. So, the theorem does not extend to arbitrary finitely generated $R$-modules (however, the equivalences (3)-(7) still hold). We will explore when it does presently. \qed
\end{remk}

We next consider two specific classes of commutative Artinian rings: those which are semisimple; and those which are $PIR$s. For the first class, we get a trivial corollary to Theorem 3.5, since all modules over such rings are injective. Hence, each module is its own injective envelope. In the second class (which contains the first) we get an intimate connection between bi-Hopfian modules and their necessarily bi-Hopfian injective envelopes, given in the second part of Theorem 3.9. 

\begin{cor}
If $R$ is semisimple Artinian and $M$ is an $R$-module, then the equivalent conditions of  Theorem 3.5 hold for $M.$
\end{cor}

\begin{proof}
The module $M$ is an injective $R$-module.
\end{proof}

\begin{remk}
For $R$ semisimple Artinian, we may add additional items to the list of equivalent conditions in Theorem 2.5: (9) $M$ is finitely generated injective; (10) $M$ is finitely generated projective; (11)  $M$ is $R$-reflexive; and (12) $M$ is linearly compact. A semisimple ring $R$ is also a {\em quasi-Frobenius} ring ($QF$), a Noetherian ring (necessarily Artinian) that is  injective as an $R$-module. Any $QF$ ring $R$ is a {\em cogenerator} for its category of modules, that is, any $R$-module $M$ can be embedded in a direct product of copies of $R.$ So, a $QF$ ring $R$ defines a Morita duality from the category of $R$-modules to itself \cite [Cor. 19.44]{Lam2}. We call an $R$-module $M$ {\em $R$-reflexive} if $M^{**} \cong M,$ where $M^*=\mathrm{Hom}_R(M,R)$ is the {\em $R$-dual} of $M.$ By \cite[Thm. 24.8]{A-F}, $M$ is $R$-reflexive if and only if $M$ is finitely generated if and only if $M$ is finitely cogenerated.  

There is a nice discussion of linearly compact modules in \cite[\S 19F, p 527]{Lam2} which includes an algebraic definition of linearly compact in terms of solutions of systems of congruences. According to Lam, the connection between reflexivity and linear compactness was not noticed until M{\"{u}}ller's paper \cite {Mull}. The relevant fact in our context  is that if $R$ is  Artinian, then an $R$-module $M$ is linearly compact if and only if $M$ has DCC if and only if  $M$ is finitely generated.

In any event, for commutative semisimple $R$ the equivalences we have mentioned subsume a result of Hiremath \cite[Theorem 16]{Hiremath} that if $R$ is semisimple Artinian, then an $R$-module $M$ is Hopfian if and only if it has finite length.

Nakayama introduced $QF$ rings in 1938 \cite[p 333]{A-F} (or 1939 \cite[p 287]{A-F}; see also footnote 1 of \cite{Nak1} and the introduction to \cite{Nak2}). Dieudonn{\'{e}} \cite{Jean} studied what he called ``perfect" duality, building off the basic theory of duality set forth by Bourbaki \cite[\S II.7.5]{Bour}. He added four properties, all satisfied by finite dimensional vector spaces, to constitute this perfect duality (one is reflexivity). Ultimately, he showed that all finitely generated $R$-modules have perfect duality if and only if $R$ is $QF.$  As semisimple rings are $QF,$ the modules of Corollary 3.7  have perfect duality. There is a discussion of the basic duality properties of $QF$ rings in \cite[\S 15C, pp 414-417]{Lam2}. For Morita duality, there is a detailed treatment in \cite [\S 19]{Lam2} and in \cite[\S 24]{A-F}.

For us, it was interesting to see that if $R$ is semisimple Artinian, the classes of Hopfian, co-Hofian, finitely generated projective, and finitely generated injective modules coincide. $\qed$
\end{remk}

It is well-known that if $R$ is an Artinian $PIR$, then an $R$-module $M$ is a direct sum of cyclic $R$-modules, each isomorphic to $R/\mathfrak{m}^k$ for some maximal ideal $\mathfrak{m}$ of $R$ and positive integer $k$ satisfying $1 \leq k \leq s_{\mathfrak{m}},$ $s_{\mathfrak{m}}$ being the least positive integer for which $\mathfrak{m}^s = \mathfrak{m}^{s+1} = \ldots $ (see, for example \cite[Thm 6.9]{S-V}).  For a given $\mathfrak{m},$ let $M_{\mathfrak{m}}=\oplus_{k=1}^{s_{\mathfrak{m}}}(\oplus_{n_k} R/\mathfrak{m}^k),$  so that $M = M_{\mathfrak{m_1}} \oplus \ldots \oplus M_{\mathfrak{m_t}},$ where the $\mathfrak{m}_i$ are the maximal ideals of $R$. It follows that $M$ is co-Hopfian (resp. Hopfian) if and only if each $M_{\mathfrak{m}}$ is (Theorem 1.1).  Hence, $M$ is bi-Hopfian if and only if each $M_{\mathfrak{m}}$ is.

If $M_{\mathfrak{m}}$ is co-Hopfian (resp. Hopfian), then each $n_k$ is finite (Theorem 1.1.3), so that $M_{\mathfrak{m}}$ is finitely generated. So, $M_{\mathfrak{m}}$ bi-Hopfian implies $M_{\mathfrak{m}}$ finitely generated. Conversely, since $R$ is Artinian, $M_{\mathfrak{m}}$ finitely generated implies that $M_{\mathfrak{m}}$  has finite length, implying that $M_{\mathfrak{m}}$ is both co-Hopfian and Hopfian, hence bi-Hopfian. So, we have proved the first part of the next result.

\begin{thm}
Let $R$ be an Artinian $PIR$ and $M$ an $R$-module.
\begin{enumerate}
\item $M$ is bi-Hopfian if and only if $M$  is finitely generated.
\item $M$ is bi-Hopfian if and only if $E(M)$ is bi-Hopfian.
\end{enumerate}
\end{thm}

\begin{proof}
We need only prove (2). Using Theorems  3.3 and 3.5, (1),  and Hopkins-Levitzki, we see that  $M$ is bi-Hopfian $\Leftrightarrow$  $M$ is finitely generated $\Leftrightarrow$ $E(M)$ is finitely generated $\Leftrightarrow$ $E(M)$ is bi-Hopfian.
\end{proof}

\begin{cor}
Let $R$ be an Artinian $PIR$ and $M$ an $R$-module with injective envelope $E=E(M).$ Any one of the equivalent conditions of Theorem 3.5 holds for $M$ if and only if it holds for $E.$ 
\end{cor}

It would be nice if part $(1)$ of Theorem 3.9 held for Artinian rings generally. Unfortunately, this is not the case. In the class of $QF$ rings mentioned earlier, there are rings which have co-Hopfian modules which are not finitely generated. 

\begin{illus}
We thank E. Enochs for walking us through some of the finer points of this example which was graciously provided by J. Rickard (we have modified Rickard's argument somewhat).

Let $k$ be an infinite field and $A=k[x,y]/(x^2,y^2).$ Note that $A$ is a four-dimensional $k$-algebra with $k$-basis $\{1,x,y,xy\}$ (for convenience, we dispense with the usual ``bar" notation for the cosets of $(x^2,y^2)$  represented by $x$ and $y).$  

In $A,$ let $\mathfrak{m}=(x,y).$ It is straightforward to show that $\mathfrak{m}^3=(0),$ and that $A/\mathfrak{m} \cong k.$ Hence, $\mathfrak{m}$ is a maximal ideal in $A.$ Since $(0) = \mathfrak{m}^3$ is contained in any prime ideal $P$ of $A,$ it follows that $P= \mathfrak{m},$ so that $\mathfrak{m}$ is the unique prime ideal in $A.$ So, $(A,\mathfrak{m})$ is local Artinian. By \cite[Example 15B, p.68, and Prop 3.14]{Lam2}, $A$ is $QF.$

 For $\lambda$ in $k^*=k-\{0\},$ the $A$-module $(x-\lambda y)A$ is a two-dimensional $k$-subalgebra of $A$ and has $k$-basis $\{x-\lambda y, xy\}.$ Since $\{1,x-\lambda y,y,xy\}$ is also a $k$-basis for $A,$ the $A$-module $M_{\lambda}=A/(x-\lambda y)A$  is a two-dimensional $k$-algebra with $k$-basis $\{1,y\}.$ Let $S_{\lambda}=Ay,$ which is a $1$-dimensional subalgebra of $M_{\lambda},$ and a cyclic $A$-module. Since $\mathfrak{m}$ is the annihilator of $S_{\lambda},$ $S_{\lambda}$ is a simple $A$-module, and hence is bi-Hopfian. 

Let $S = \oplus_{k^*}S_{\lambda}.$ Clearly, $S$ is not finitely generated. However, it is co-Hopfian. Let $f:S \rightarrow S,$ be an injective $A$-module homomorphism. We can think of $f$ as represented by a matrix $(f_{\mu \lambda})$ with $f_{\mu \lambda}:S_{\lambda} \rightarrow S_{\mu}.$

If $f_{\mu \lambda}:S_{\lambda} \rightarrow S_{\mu}$ is an isomorphism (or simply an injection), then $f$ is a linear transformation of $1$-dimensional vector spaces, so is multiplication by some nonzero scalar $d \in k.$ So, $f(y)=dy$ and $f(x)=dx.$ So, $f((\mu y)/d)=\mu y,$ and $f(x/d)=x.$ Hence, $(x/d)-(\mu y)/d$ is in the kernel of $f$ since $x - \mu y=0$ in $S_{\mu}.$ Since $f$ is an injection, $x=\mu y$ in $S_{\lambda}.$ But, $x=\lambda y$ in $S_{\lambda}.$ Therefore, $\mu = \lambda.$  Since the $S_{\lambda}$ are simple, maps $S_{\lambda} \rightarrow S_{\mu}$ with $\mu \neq \lambda$ must be zero. So, any injective endomorphism of $\oplus_{\lambda \in k^*} S_{\lambda}$ will diagonalize as a product of isomorphisms $\prod f_{\lambda \lambda}$ (since the $S_{\lambda}$ are simple) and so be an isomorphism of $S.$         \qed
\end{illus}

\section{The case of an arbitrary Artinian PIR.}

Since Hopkins-Levitzki is valid for arbitrary Artinian rings $R$, it is natural to ask to what extent the results of section 3 can be duplicated if the ring $R$ is not assumed commutative. There is one potentially complicating feature in the noncommutative case, mentioned in the introduction to section 3: indecomposable injectives need not be finitely generated. However, if $R$ is an Artinian $PIR$, then this problem vanishes since $R$ will be $QF.$ This fact allows us to generalize the results for injective modules in section 3 to the noncommutative setting, at least when $R$ is an Artinian $PIR.$

Henceforth, unless otherwise specified, we assume that $R$ is a left Artinian ring, and that all modules are left modules. We will consider two specific classes of Artinian rings: those which are semisimple; and those which are $PIR$s (the second class contains the first). If $R$ is semisimple, we get a trivial corollary to Theorem 3.5, since all modules (left or right) over such rings are injective, with the simple modules being the indecomposable injectives. Since there are only finitely many isomorphism classes of simple modules, we can then appeal to Theorem 1.1 and \cite[Thm 3.13 Cor]{S-V}, the latter a useful corollary to the Krull-Remak-Schmidt-Azumaya Theorem (part of the machinery used by Matlis to deduce his structure result for injectives over a commutative Noetherian $R$).

\begin{thm}
If $R$ is semisimple Artinian and $M$ is an $R$-module, then the equivalent conditions of  Theorem 3.5 hold for $M.$ \qed
\end{thm}

Before proceeding, we point out the following result, which gives a sufficient condition for a module to be Hopfian (co-Hopfian).

\begin{thm}
Let $R$ be an Artinian ring with radical $J,$ and let $M$ be an $R$-module. If $M/JM$ is Hopfian (co-Hopfian), then so is $M.$
\end{thm}

\begin{proof}
The module $M/JM$ is an $R/J$-module, and $R/J$ is semisimple since $R$ is Artinian \cite[Prop 15.17]{A-F}. But then $M/JM$ is semisimple \cite[Cor 15.18]{A-F}. Let $S_1, \ldots, S_m$ be a complete set of representatives of the isomorphism classes of simple $R/J$-modules. Then
$$M/JM \cong \oplus_{i=1}^{m}(\oplus_{k_i}S_i).$$
By Theorem 1.1, $M/JM$ Hopfian (co-Hopfian) as an $R/J$-module if and only if each $k_i$ is finite if and only if $M/JM$ is finitely generated. But $R/J$ is Artinian since $R$ is. Therefore, $M$ is finitely generated \cite[Cor 15.21]{A-F}. Hence, $M$ is Hopfian (co-Hopfian) as an $R$-module since $R$ is Artinian (and hence Noetherian) .
\end{proof}

For the case of an Artinian $PIR$  (a ring that is both left and right Artinian and for which each left or right ideal is principal), the main ideas supporting our argument may be found in sections 31, 32, 25, 27, and 28 of \cite{A-F}, and in results from earlier sections mentioned in these five.  

If $R$ is any Artinian ring with Jacobson radical $J,$ then $J$ is nilpotent, and $R/J$ is semisimple. Hence, $R$ is a {\em semiprimary} ring. Moreover, idempotents of $R/J$ lift to idempotents of $R,$ so that $R$ is also {\em semiperfect}. This last fact allows us to characterize the indecomposable projective $R$-modules and the simple $R$-modules.

We recall quickly a few facts about idempotents that are scattered throughout \cite{A-F}. In any ring $R,$ an element $e$ is {\em idempotent} if $e^2=e.$ Idempotents $e_1,e_2$ (usually nonzero) are {\em orthogonal} if $e_1e_2=0=e_2e_1.$ Note that if $e$ is idempotent, then so is $1-e$ and the two are orthogonal. A nonzero idempotent $e$ is {\em primitive} if it is not the sum of nonzero orthogonal idempotents.

For $R$ semiperfect, an $R$-module $M$ is {\em primitive} if $M \cong Re$ for some primitive idempotent of $R.$  A set $e_1, \ldots, e_m$ of idempotents of $R$ is {\em basic} in case the $e_i$ are pairwise orthogonal and  the modules $Re_1, \ldots, Re_m$ form a complete and irredundant set of repesentatives of the isomorphism classes of primitive $R$-modules (so that the $e_i$ are primitive idempotents). Every complete set of pairwise orthogonal primitive idempotents contains a basic set and if $e_1, \ldots, e_m$ is such a basic set, then $Re_1, \ldots, Re_m$ is a complete and irredundant set of representatives for the isomorphism classes of projective  $R$-modules, and $Re_1/Je_1, \ldots, Re_m/Je_m$ is a complete and irredundant set of representatives for the isomorphism classes of simple $R$-modules (\cite[Prop 27.10]{A-F}). So, we have a structure theorem \cite[Thm 27.11]{A-F} for projectives over a semiperfect ring that resembles Matlis' theorem for injectives over a commutative Noetherian ring (Theorem 2.1).

\begin{thm}
Let $R$ be semiperfect, and $e_1, \ldots , e_m$ a basic set of primitive idempotents of $R.$ If $P$ is a projective $R$-module, then there are unique cardinal numbers $n_1, \dots , n_m$ for which $$P \cong \oplus_{i=1}^m (\oplus_{n_i}Re_i).$$ The cardinals $n_i$ characterize $P$ up to isomorphism and some are possibly $0.$
\end{thm}

If $R$ is Artinian, then any $R$-module $M$ is an essential extension of its socle, equivalently, every nonzero submodule of $M$ contains a simple submodule \cite[Prop 3.17 and Cor]{S-V}. Hence, the indecomposable injective $R$-modules are 
$$E(Re_1/Je_1), \ldots, E(Re_m/Je_m).$$ 

Eisenbud and Griffith \cite[Prop 2.3]{E-G2} showed that if $R$ is an Artinian $PIR,$ then $R$ is both a serial ring and $QF.$ So, the latter condition implies that an $R$-module $M$ is projective if and only if it is injective. So, the indecomposable injectives $E(Re_i/Je_i)$ are the $Re_i$ is some order. Consequently, we have the injective analogue of Theorem 2.1 and a characterization of those injectives that are Hopfian (co-Hopfian).

\begin{thm}
Let $R$ be an Artinian $PIR,$  $e_1, \ldots , e_m$ a basic set of primitive idempotents of $R,$ and $E$ an injective $R$-module. Then then there are unique cardinal numbers $n_1, \ldots , n_m$ for which $$E=\oplus_{i=1}^m (\oplus_{n_i}Re_i),$$ and the cardinal numbers $n_i$ characterize $E$ up to isomorphism. Moreover, $E$ is Hopfian (co-Hofian) if and only if $n_i$ is finite for all $i$ (some $n_i$ may be $0$). Hence, $E$ is Hopfian if and only if it is co-Hopfian.
\end{thm}

\begin{proof}
If all $n_i$ are finite $E$ is finitely generated, so $E$ is both Hopfian and co-Hopfian. Conversely, a now familiar argument shows that if $R$ is Hopfian (co-Hopfian), then each $n_i$ must be finite. The final statement is obvious.
\end{proof}

\begin{cor}
An injective module over an Artinian $PIR$ is Hopfian (co-Hopfian) if and only if it is finitely generated. \qed
\end{cor}

The theorem corresponding to Theorem 3.5 is the following.

\begin{thm}
Let $E$ be an injective $R$-module, $R$ an Artinian $PIR.$ The equivalent conditions of Theorem 3.5 hold for $E.$ 
\end{thm}

\begin{proof}
Recall that each $Re_i$ is one of the $E(Re_j/Je_j).$
\end{proof}

\begin{remk}
Since an Artinian $PIR$ is $QF,$ the conditions $(8)-(12)$ added in Remarks 3.6, 3.8 are also valid for $E.$
\end{remk}

It is clear from the corollary that if $M$ is an $R$-module and $E(M)$ is Hopfian (co-Hopfian), then $M$ is finitely generated and so is Hopfian (co-Hopfian). Thus, $E(M)$ finitely generated implies that $M$ is finitely generated. The converse is also true. Since $R$ is Artinian it is Noetherian, so if $M$ is a finitely generated $R$-module, then $E(M)$ is a finite direct sum of indecomposable injective $R$-modules \cite[Cor 3.50]{Lam2}. But $R$ is also a $PIR$, hence $QF,$  so, as shown above, the indecomposable injectives are the $Re_i,$ and these are cyclic. Therefore, we have the following analogue of Theorem 3.3.

\begin{thm}
Let $R$ be an Artinian $PIR$ and $M$ an $R$-module. Then $M$ is finitely generated if and only if $E(M)$ is finitely generated.
\end{thm}

If we can show that $M$ is Hopfian (co-Hopfian) if and only if $M$ is finitely generated, then we will be able to conclude that $M$ is Hopfian (co-Hopfian) if and only if $E(M)$ is Hopfian (co-Hopfian), and Corollary 3.10 will be valid  for modules over an arbitrary Artinian $PIR.$

The previously mentioned result of Eisenbud-Griffith established that if $R$ is an Artinian $PIR,$  then $R$ is a serial ring. Any module over a serial ring is a direct sum of uniserial modules. However, in \cite[Thm 2.1]{E-G2}, they also established that the uniserial modules may be assumed {\em homogeneous} (a uniserial module is homogeneous if its composition factors are mutually isomorphic). Any Artinian ring has only finitely many isomorphism classes of simple modules. For $R$  an Artinian  $PIR$ they will be the $Re_i/Je_i$ mentioned earlier. For simplicity, call these simple modules $S_1, \ldots, S_m.$ Let the $R$-module $M$ have decomposition $$M=\oplus_{i \in I}M_i$$ with each $M_i$ homogeneous uniserial. 

To mimic the proof of Theorem 3.9.1, we would begin by letting $H_j$ represent the sum of those $M_i$ whose simple submodule is  isomorphic to $S_j.$ We could continue if we knew that $R$ has only finitely many isomorphism classes of uniserial module, that is, $R$ has {\em finite uniserial type}. As far as we have been able to discover, determining if $R$ has finite uniserial type is an open problem. It is even an open problem at the level of Artin algebras, algebras over a commutative Artinian ring $S$ that are finitely generated as $S$-modules (see Problem 2, p. 411, in the list of open problems listed in \cite{ARS}). D'este, Kaynarca, and T\"{u}t\"{u}nc\"{u} \cite{D'e} give a summary of progress in the Artin algebra case in their introduction, and present an example of an Artin algebra having two nonisomorphic uniserial modules of length two having the same composition factors. B. Huisgen-Zimmermann \cite{BH-Z} provides the tools to answer the question in the case in which the ring $S$ is an algebraically closed field. So, for the moment at least, the rest of the results of Section 3 remain unresolved for general Artinian $PIR$s.

\section{Epilogue}

We close with brief remarks for those who might be curious about the relevance of Artinian $PIR$s and to mention some of the literature  pertaining to the results we used from \cite{A-F}.

We have a classical result that a {\em bounded} abelian group $G$ is a direct sum of cyclic groups (attributed by Fuchs \cite[Thm 17.2]{fuchs} to Pr\"{u}fer \cite{Prufer} and Baer \cite{Baer}). Here, bounded means the group has bounded exponent. That is, there is a least positive integer $n$ such that $nx=0$ for all $x \in G.$ A bounded group is naturally a module over the ring $\mathbb{Z}/n\mathbb{Z},$ which is an Artinian $PIR,$ uniserial if $n$ is a prime power  and serial otherwise. A cyclic group of order $n$ behaves similarly, it is uniserial if $n$ is a prime power and serial otherwise.

 In his influential monograph on infinite abelian groups,   Kaplansky \cite[p78]{Kap} recounts the work of K\"{o}the \cite{Kothe}, Asano \cite{Asano},  Cohen and Kaplansky \cite{C-K}, Chase \cite{Chase}, and Faith and Walker \cite{F-W} extending the direct sum of cyclics result to principal ideal rings. There is a nice discussion of the commutative Artinian $PIR$ case in \cite[Section 6.1]{S-V}. So, Artinian $PIR$s provide a natural class of rings in which to try to generalize certain results on abelian groups.

Nakayama \cite[Thm 17]{Nak2} showed that if $R$ is left Artinian, then $R$ is serial if and only if every finitely generated left $R$-module is a direct sum of uniserial modules. Eisenbud and Griffith \cite[Thm 1.2]{E-G1} removed the finitely generated restriction. In the subsequent paper \cite[Thm 2.1]{E-G2}, they went on to show  that $R$ is an Artinian $PIR$ if and only if every left $R$-module is a direct sum of homogeneous uniserial modules.  Later in that same paper (Proposition 2.3), they showed that if $R$ is an Artinian $PIR,$ then $R$ is serial and $QF.$

Finally, Faith and Walker \cite[Thms 5.3, 5.5]{F-W} proved two results that were important to us: A ring is $QF$ if and only if every injective module is projective if and only if every projective module is injective; and, a ring is $QF$ if and only if each injective right module is a direct sum of cyclic modules that are isomorphic to principal indecomposable right ideals (Theorem 4.4 is the left-hand version).

\section{Acknowledgments}

I wish to thank three people who have had significant influence on me and on my mathematics: E. E.  Enochs, who taught me graduate algebra and who sparked my interest in module theory and abelian groups, he remains a source of counsel; the late E. D. Davis, who taught me  more algebra, and introduced me to commutative algebra; and my thesis advisor, H. I. Brown, who kindled in me an affection for summability theory that lingers to this day. Finally, I give special thanks to Margaret Carney OSF for her support at a critical stage of my career.

\end{document}